\documentclass[runningheads]{cl2emult}

\usepackage{makeidx}  
\usepackage{graphicx} 
\usepackage{subeqnar} 
\usepackage{multicol} 
\usepackage{cropmark} 
\usepackage{math}     
\makeindex            

\begin{document}

\title*{Umbral nature of the Poisson random variables}
\titlerunning{Umbral nature of the Poisson random variables}

\author{Elvira Di Nardo
\and Domenico Senato}
\authorrunning{E. Di Nardo, D. Senato}
\institute{Dipartimento di Matematica, Universit\`a degli Studi
della Basilicata, \\
Via N. Sauro 85, 85100 Potenza, Italy\\
\tt{\{dinardo,\,senato\}@unibas.it}}
\maketitle
\begin{abstract}
Extending the rigorous presentation of the \lq \lq classical umbral
calculus\rq\rq$\,$ \cite{RotaTaylor}, the so-called partition
polynomials are interpreted with the aim to point out the
umbral nature of the Poisson random variables. Among the new umbrae
introduced, the main tool is the partition umbra that leads also to a simple
expression of the functional composition of the exponential power
series. Moreover a new short proof of the Lagrange inversion formula is
given.
\end{abstract}
%
\section{Introduction}
The symbolic method, nowadays known as {\it umbral calculus},
has been extensively used since the nineteenth century although
the mathematical community was sceptic of it, maybe owing to its
lack of foundation. This method was fully developed by Rev. John Blissard
in a series of papers beginning from 1861 \cite{Blissard1}$\div$\cite{Blissard11},
nevertheless it is impossible to attribute
the credit of the originary idea just to him since the
Blissard's calculus has a mathematical source in the symbolic
differentiation. In \cite{Lucas} Lucas even claimed that
the umbral calculus has its historical roots in the writing
of Leibniz for the successive derivatives of a product with
two or several factors. Moreover Lucas held that this symbolic method
had been subsequently developed by Laplace, by Vandermonde, by Herschel and
augmented by the works of Cayley and of Sylvester in the
theory of forms. Lucas's papers attracted considerable
attention and the predominant contribution of Blissard to
this method was kept in the background. Bell reviewed
the whole subject in several papers, restoring the purport
of the Blissard's idea \cite{Bell4} and in 1940 he
tried to give a rigorous foundation of
the mystery at the ground of the umbral calculus \cite{Bell2}. It was
Gian-Carlo Rota \cite{Bellnumbers} who twenty-five years
later disclosed the \lq \lq umbral magic art\rq\rq$\,$ of lowering
and raising exponents bringing to the light the underlying
linear functional. In \cite{Mullin} and \cite{VIII} the ideas
from \cite{Bellnumbers} led Rota and his collaborators to
conceive a beautiful theory originating a large variety of
applications. Some years later, Roman and Rota gave
rigorous form to the umbral tricks in the setting of the
Hopf algebra. On the other hand, as Rota himself has written \cite{RotaTaylor}:
\lq \lq ...Although the notation of Hopf algebra satisfied the
most ardent advocate of spic-and-span rigor, the translation of
\lq \lq classical\rq\rq$\,$ umbral calculus into the newly found
rigorous language made the method altogether unwieldy and
unmanageable. Not only was the eerie feeling of witchcraft lost
in the translation, but, after such a translation, the use of
calculus to simplify computation and sharpen our intuition
was lost by the wayside...\rq\rq$\,$ Thus in 1994 Rota and Taylor
\cite{RotaTaylor} started a rigorous and simple presentation of the umbral
calculus in the spirit of the founders. The present article
refers to this last point of view.
\par
As it sometimes happens in the practice of the mathematical
investigation, the subject we deal with does not develop
the originary idea from which our research started in
the spring of 1997, but this paper is closely related
to it. In that period, Gian-Carlo Rota was
visiting professor at the University of Basilicata and, during
one of our latest conversations before his leaving, he shared
with us his close interest for a research project: a combinatorial
random variable theory. The delicate question arising from the
underlying foundation side and the left short time led us to
protract the discussion via e-mail intertwining it with
different activities for several months. The following year,
Gian-Carlo Rota held his last course in Cortona and we did not
miss the opportunity to spend some time with him. We resumed
the thread of our conversations and presented him with the doubts
that gradually took hold of us. As usually, his contribution
disclosed new horizons that have led us to write these
pages.
\par
Our starting point is the umbral notion of the Bell numbers.
Many classical identities relating to these numbers are
expressed in umbral notation attaining up to a new umbra,
the partition umbra, connected with the so-called \lq \lq partition
polynomials\rq\rq$\,$ generated by expanding the exponential function $\exp(f(x))$ into
an exponential power series. The hereafter developed theory of the
Bell umbrae is not only an example of the computational power of the
umbral calculus but it offers, we would like to believe, a natural
way to interpret the functional composition of exponential
power series tested by a new proof of the Lagrange inversion formula.
Here the point operations extended with a new one play
a central role. From a probabilistic point of view, the
functional composition of exponential power series is closely
related to the family of Poisson random variables so
that these random variables have found a natural umbral
interpretation through the Bell umbrae. In particular the
probabilistic counterpoint of the partition umbra is the
compound Poisson random variable. Also the less familiar
randomized Poisson random variable unexpectedly find an
umbral corresponding and in turn the umbral composition
gives the way to generalize this last Poisson random variable.
\par
What's more, we believe that the probabilistic interpretation
of the partition umbra and of the umbral composition could
give a probabilistic meaning to the Joyal species theory
\cite{Joyal}, namely a combinatorial random variable theory
that we hope to deal in forthcoming publication.
%
\section{The classical umbral calculus}
We take a step forward in the program of the rigorous foundation
of the classical umbral calculus initiated by Rota and Taylor
\cite{SIAM}, \cite{RotaTaylor}, \cite{Taylor}. \par
In the following we denote by $R$ a commutative integral domain
whose quotient field is of characteristic zero and
by $A=\{\alpha,\beta, \ldots \}$ a set whose elements are called {\em umbrae}.
An umbral calculus is given when is assigned a linear functional
$E: R[A,x,y] \rightarrow R[x,y]$ such that:
\begin{description}
\item[{\it i)}] $E[1]=1;$
\item[{\it ii)}] $E[\alpha^i \beta^j \cdots \gamma^k x^n y^m] =
x^n y^m E[\alpha^i]E[\beta^j] \cdots E[\gamma^k]$
for any set of distinct umbrae in $A$ and for $i,j,\ldots,k,n,m$
nonnegative integers (uncorrelation property);
\item[{\it iii)}] it exists an element $\epsilon \in A$  such that
$E[\epsilon^n] = \delta_{0,n},$ for any nonnegative integer $n,$
where
$$\delta_{i,j} = \left\{ \begin{array}{cc}
1 & \hbox{if $i=j$} \\
0 & \hbox{if $i \ne j$}
\end{array} \right. \,\, i,j \in N;$$
\item[{\it iv)}] it exists an element $u \in A$  such that $E[u^n]=1,$
for any nonnegative integer $n.$
\end{description}
The umbra $\epsilon$ is named {\it augmentation} as Roman and Rota first called
it \cite{umbral}. We will call the umbra $u$ the {\it unity} umbra.
\par
A sequence $a_0,a_1,a_2, \ldots$ in $R[x,y]$ is said to be
umbrally represented by an umbra $\alpha$ when
$$E[\alpha^i]=a_i, \quad \hbox{for} \,\, i=0,1,2,\ldots$$
so that the linear functional $E$ plays the role of an evaluation
map.
\par
As Rota suggested, there is an analogy between umbrae and random
variables (r.v.) (see \cite{Taylor}), so we will refer to the
elements $a_i$ in $R[x,y]$ as {\it moments} of the umbra $\alpha.$
The umbra $\epsilon$ can be view as the r.v. which takes the value
$0$ with probability 1 and the umbra $u$ as the r.v. which
takes the value $1$ with probability 1. \par
An umbra is said to be a scalar umbra if the moments are elements of
$R$ while it is said to be a polynomial umbra if the moments are
polynomials. Note that if the sequence $a_0, a_1, a_2, \ldots$
is umbrally represented by a scalar umbra $\alpha,$ then it is $a_0=1.$
In the same way, for polynomial umbrae, a sequence
of polynomials $p_0(x),p_1(x),p_2(x), \ldots$ will always
denote a sequence of polynomials with coefficients in $R$ such
that $p_0(x)=1$ and $p_n(x)$ is of degree $n$ for every
positive integer $n.$
\par
A polynomial $p \in R[A]$ is called an {\it umbral polynomial}. The {\it support} of
$p$ is defined to be the set of all occurring umbrae of $A$.
Two umbral polynomials are said to be {\it uncorrelated} when their support are
disjoint.
\par
If $\alpha$ and $\beta$ are either scalar either polynomial
umbrae, we will say that $\alpha$ and $\beta$ are {\it umbrally
equivalent} when
$$E[\alpha]=E[\beta],$$
in symbols $\alpha \simeq \beta.$ Two scalar (or polynomial) umbrae are said
to be {\it similar} when $E[\alpha^k]=E[\beta^k], \, k=0,1,2,\ldots$ or
$$\alpha^k \simeq \beta^k, \, k=0,1,2,\ldots$$
in symbols $\alpha \equiv \beta.$ The notion of equivalence and similarity
for umbral polynomials is obvious.
\par
The formal power series
$$ e^{\alpha t}=u + \sum_{n \geq 1} \alpha^n \frac{t^n}{n!}$$
is said to be the {\it generating function} of the umbra
$\alpha.$ Moreover, if the sequence $a_0,a_1,a_2,\ldots$ has
(exponential) generating function $f(t)$ and is umbrally represented
by an umbra $\alpha$ then $E[e^{\alpha t}]=f(t),$ in symbols
$e^{\alpha t}\simeq f(t).$ When $\alpha$ is regarded as a r.v.,
$f(t)$ is the moment generating function. The notion of equivalence
and similarity are extended coefficientwise to the generating functions
of umbrae so that $\alpha \equiv \beta$ if and only if $e^{\alpha t}
\simeq e^{\beta t}.$  Note that $e^{\epsilon t} \simeq 1$ and $e^{u t}
\simeq e^x.$
%
\subsection{The point operations}
%
The notion of similarity among umbrae comes in handy in order to
express sequences such
\begin{equation}
\sum_{i=0}^n \left( \begin{array}{c} n \\ i \end{array} \right)
a_{i}a_{n-i}, \,\, n=0,1,2,\ldots
\label{(eq:1)}
\end{equation}
as moments of umbrae. The sequence (\ref{(eq:1)}) cannot be represented
by using only the umbra $\alpha$ with moments $a_0,a_1,a_2, \ldots$ because
$a_i a_{n-i}$ could not be written as $E[\alpha^i \alpha^{n-i}],$
being $\alpha$ related to itself.
If we will assume that the umbral calculus we deal is
saturated \cite{SIAM}, the sequence $a_0,a_1,a_2,\ldots$ in
$R[x,y]$ is represented by infinitely many distinct (and thus similar)
umbrae. Therefore, if we choose two similar umbrae $\alpha, \alpha^{\prime}$,
they are uncorrelated and
$$\sum_{i=0}^n \left( \begin{array}{c} n \\ i \end{array} \right)
a_{i}a_{n-i} = E\left[\sum_{i=0}^n \left( \begin{array}{c} n \\ i \end{array} \right)
\alpha^i \alpha^{n-i}\right]= E[(\alpha + \alpha^{\prime})^n].$$
Then the sequence (\ref{(eq:1)}) represents the moments of the
umbra $(\alpha+\alpha^{\prime}).$ This matter was first
explicitly pointed out by E.T. Bell \cite{Bell2} who was not
able to provide an effective notation:
\par \bigskip
\footnote{The quotation needs more details.
It is
$$\alpha x \dot{+} \cdots \dot{+} \varepsilon x \equiv (\alpha x_0 + \cdots +
\varepsilon x_0, \ldots, \alpha x_N + \cdots + \varepsilon x_N, \ldots),$$
where $x$ is an umbra. The formula (1.22) is
$$(\alpha a \dot{+} \cdots \dot{+} \varepsilon x)^N =
\sum M_{S_1, \ldots, S_T} \alpha^{S_1} \cdots \varepsilon^{S_T}
a^{S_1} \cdots x^{S_T}$$
with $M_{S_1, \ldots, S_T}$ the coefficient of
$x^{S_1}_1 \cdots x^{S_T}_T$ in the expansion of $(x_1+\cdots+x_T)^N$
through the multinomial theorem. The formula (1.20) is
$$(\alpha a \dot{+} \cdots \dot{+} \varepsilon x)^N \equiv
\sum M_{S_1, \ldots, S_T} \alpha^{S_1} \cdots \varepsilon^{S_T}
a_{S_1} \cdots x_{S_T}.$$}{\footnotesize \lq\lq...If in
$\alpha x \dot{+} \cdots \dot{+} \xi x$ there are precisely $T$
summands $\alpha x, \ldots, \xi x$ each of which is a scalar product
of a scalar and $x$, we replace $(\rightarrow)$ the $T$ $x'$s by $T$ distinct
umbrae, say $a,\ldots, x$ in any order, and indicate this
replacement by writing
$$\alpha x \dot{+} \cdots \dot{+} \xi x \rightarrow
\alpha a \dot{+} \cdots \dot{+} \xi x.$$
Then $(\alpha a \dot{+} \cdots \dot{+} \xi x)^N$ is to be calculated by
(1.22) and the exponents are degraded as in (1.20).
In the result, each of $a, \ldots, x$ is replaced $(\leftarrow)$
by $x;$ the resulting polynomial is defined to be the $N-$th
power $(\alpha x \dot{+} \cdots \dot{+} \xi x)^N$ of
the sum $(\alpha x \dot{+} \cdots \dot{+} \xi x).$
For example
\begin{eqnarray*}
(\alpha x + \beta x)^3 & \rightarrow & (\alpha a + \beta x)^3;
\\
(\alpha a + \beta x)^3 & = & \alpha^3 a_3 x_0 + 3 \alpha^2
\beta a_2 x_1 + 3 \alpha \beta^2 a_1 x_2 + \beta^3 a_0 x_3, \\
& \leftarrow & \alpha^3 x_3 x_0 + 3 \alpha^2 \beta x_2 x_1 + 3
\alpha \beta^2 x_1 x_2 + \beta^3 x_0 x_3; \\
(\alpha x + \beta x)^3 & = & (\alpha^3+\beta^3) x_0 x_3 + 3
\alpha \beta (\alpha + \beta) x_1 x_2. \hbox{...\rq\rq}
\end{eqnarray*}}
\par
The last identity makes sense when the left side is replaced
by $(\alpha x + \beta x^{\prime})^3$ with $x \equiv x^{\prime},$
but Bell did not have the notion of similar umbrae. However, the need
of handling sequences like (\ref{(eq:1)}) leads to introduce
some new operations between umbrae, as showed in the next
sections.
%
\subsection{The point product}
%
We shall denote by the symbol $n.\alpha$ an auxiliary umbra similar to
the sum $\alpha^{\prime}+\alpha^{\prime \prime}+ \ldots + \alpha^{\prime\prime\prime}$
where $\alpha^{\prime},\alpha^{\prime \prime},\ldots,\alpha^{\prime\prime
\prime}$ are a set of $n$ distinct umbrae each of which is similar to the
umbra $\alpha.$ We assume that $0.\alpha$ is an umbra similar
to the augmentation $\epsilon.$ A similar notion $n.p$ is
introduced for any umbral polynomial. The following statements
are easily to be proved:
\begin{proposition}
\begin{description}
\item[\it{(i)}] If $n.\alpha \equiv n.\beta$ for some integer $n \ne 0$
then $\alpha \equiv \beta;$
\item[\it{(ii)}] if $c \in R$ then $n.(c \alpha) \equiv c (n.\alpha)$
for any nonnegative integer $n;$
\item[\it{(iii)}] $n.(m.\alpha) \equiv (nm).\alpha \equiv m.(n.\alpha)$
for any two nonnegative integers $n,m;$
\item[\it{(iv)}] $(n+m).\alpha \equiv n.\alpha + m.\alpha^{\prime}$
for any two nonnegative integers $n,m$ and any two distinct umbrae $\alpha \equiv \alpha^{\prime};$
\item[\it{(v)}] $(n.\alpha + n.\beta)\equiv n.(\alpha+\beta)$
for any nonnegative integer $n$ and any two distinct umbrae $\alpha$ and $\beta.$
\end{description}
\end{proposition}
\begin{proposition}
If $\alpha$ is an umbra with generating function
$e^{\alpha t} \simeq f(t),$ then the umbra $n.\alpha$ has generating function
$e^{(n.\alpha)t} \simeq  [f(t)]^n.$
\end{proposition}
\begin{proof}
It follows from the definition of the auxiliary umbra $n.\alpha.$
\qed
\end{proof}
The moments of the umbra $n.\alpha$ are the following polynomials
in the variable $n$
\begin{equation}
E[(n.\alpha)^k] = q_k(n) = \sum_{i=0}^k (n)_i B_{k,i}, \quad k=0,1,2,...
\label{(eq:2)}
\end{equation}
where $B_{k,i}=B_{k,i}(a_1,a_2,\ldots, a_{k-i+1})$ for $i\leq k$
are the (partial) Bell exponential polynomials \cite{Bell1}, $(n)_i$ is the
lower factorial and $a_i$ are the moments of the umbra $\alpha.$
Recalling that
\begin{equation}
\sum_{k=i}^{\infty} B_{k,i} \frac{t^k}{k!} = \frac{1}{i!}
[f(t)-1]^i,
\label{(eq:3)}
\end{equation}
the identity (\ref{(eq:2)}) follows from
\begin{equation}
[f(t)]^n=\sum_{i=0}^{\infty} (n)_i \frac{[f(t)-1]^i}{i!} = \sum_{k=0}^{\infty}
\left( \sum_{i=0}^k (n)_i B_{k,i} \right) \frac{t^k}{k!}.
\label{(eq:4)}
\end{equation}
If in (\ref{(eq:2)}) set $\alpha = u,$ then $q_k(n)=n^k.$
Note that $q_0(n)=1,$ $q_k(0)=0$ and the
polynomial sequence $\left\{ q_k(n) \right\}$ is of binomial
type as it follows by using the statement {\it (iv)} of
Proposition 1:
$$[(n+m).\alpha]^k \simeq [n.\alpha+m.\alpha^{\prime}]^k
\simeq \sum_{i=0}^k \left( \begin{array}{c}
k \\
i
\end{array}\right) (n.\alpha)^i (m.\alpha^{\prime})^{k-i}.$$
Moreover a variety of combinatorial identities could be umbrally
interpreted. As instance in point, the classical Abel identity
becomes
\begin{equation}
(\alpha+\beta)^n \simeq \sum_{k=0}^n \left( \begin{array}{c}
n \\
k
\end{array}\right) \alpha (\alpha - k.\gamma)^{k-1}
(\beta+k.\gamma)^{n-k}, \quad n=0,1,2,\ldots.
\label{(Abel)}
\end{equation}
\noindent \indent
The expression of the polynomial sequence
$\{q_k(n)\}$ in (\ref{(eq:2)}) suggests a way to define the
auxiliary umbra $x.\alpha$ when $x \in R,$ however it is
impossibile to give an intrinsic definition.
Up to similarity, the umbra $x.\alpha$ is the polynomial umbra with moments
\begin{equation}
E[(x.\alpha)^k]=q_k(x)= \sum_{i=0}^k (x)_i B_{k,i} \quad k=0,1,2,....
\label{(eq:5)}
\end{equation}
Note that $q_k(x)=x^k$ when $\alpha = u.$
\begin{proposition}
If $\alpha$ is an umbra with generating function $e^{\alpha t} \simeq f(t),$
then the umbra $x.\alpha$ has generating function $e^{(x.\alpha)t} \simeq  [f(t)]^x.$
\end{proposition}
\begin{proof}
It follows from (\ref{(eq:4)}) and (\ref{(eq:3)}) with $n$ replaced by
$x.$ \qed
\end{proof}
\begin{corollary}
\begin{description}
\item[\it{(i)}] If $x.\alpha \equiv x.\beta$ for  $x \in R-\{0\}$
then $\alpha \equiv \beta;$
\item[\it{(ii)}] if $c \in R$ then $x.(c \alpha) \equiv c (x.\alpha)$
for any $x \in R;$
\item[\it{(iii)}] $x.(y.\alpha) \equiv (x y).\alpha \equiv y.(x.\alpha)$
for any $x,y \in R;$
\item[\it{(iv)}] $(x + y).\alpha \equiv x.\alpha + y.\alpha^{\prime}$
for any $x,y \in R$ and any two distinct umbrae $\alpha
\equiv \alpha^{\prime};$
\item[\it{(v)}] $(x.\alpha + x.\beta)\equiv x.(\alpha+\beta)$
for any $x \in R$ and any two distinct umbrae $\alpha$ and $\beta.$
\end{description}
\end{corollary}
\begin{theorem}
Up to similarity, each polynomial sequence of binomial type is
umbrally represented by an auxiliary umbra $x.\alpha$ and
viceversa.
\end{theorem}
\begin{proof}
>From the statement{\it (iv)} of the corollary 1, it follows that the
polynomial sequence $\{q_k(x)\}$ is of binomial type.
Viceversa, first observe that from (\ref{(eq:5)}) it is
\begin{equation}
D_x[q_k(x)]_{x=0}=a_k+F(a_1,a_2,\ldots,a_{k-1})
\label{(eq:6)}
\end{equation}
where $F$ is a function of the moments $a_1,a_2,\ldots,a_{k-1}.$
Let $\{p_k(x)\}$ be a polynomial sequence of binomial type.
Through (\ref{(eq:6)}), the moments of the umbra $\alpha$ are
uniquely determined by the knowledge of the first derivative respect to $x$ of
$p_k(x)$ evaluated in $0.$ Moreover, the sequences of first derivate
respect to $x$ of $p_k(x)$ evaluated in $0$ uniquely determines
a sequence of binomial type. \qed
\end{proof}
\noindent \indent
Similarly with it has been done for the auxiliary umbra
$x.\alpha,$ we define a point product among umbrae.
Up to similarity, the umbra $\beta.\alpha$ is an auxiliary
umbra whose moments are umbrally expressed through the umbral
polynomials $q_{\alpha,k}(\beta):$
\begin{equation}
(\beta.\alpha)^k  \simeq q_{{\alpha},k}(\beta) =
\sum_{i=0}^k (\beta)_i B_{k,i}  \quad k=0,1,2,....
\label{(eq:7)}
\end{equation}
If $\alpha$ is an umbra with generating function $e^{\alpha t} \simeq f(t),$
then the identity (\ref{(eq:4)}) could be rewritten as
\begin{equation}
[f(t)]^{\beta} \simeq \sum_{i=0}^{\infty} (\beta)_i \frac{[f(t)-1]^i}{i!} \simeq
\sum_{k=0}^{\infty} \left( \sum_{i=0}^k (\beta)_i B_{k,i} \right) \frac{t^k}{k!}
\label{(eq:8)}
\end{equation}
so that $e^{(\beta.\alpha)t} \simeq  [f(t)]^{\beta}.$ Moreover
if $\beta$ is an umbra with generating function $e^{\beta t} \simeq
g(t),$ then
$$ [f(t)]^{\beta} \simeq e^{\beta log f(t)} \simeq g\left[ \log f(t) \right].$$
This proves  the following proposition.
\begin{proposition}
If $\alpha$ is an umbra with generating function $e^{\alpha t} \simeq f(t)$
and $\beta$ is an umbra with generating function $e^{\beta t} \simeq
g(t),$  then the umbra $\beta.\alpha$ has generating function
\begin{equation}
e^{(\beta.\alpha)t} \simeq  [f(t)]^{\beta} \simeq g\left[ \log f(t) \right].
\label{(eq:9)}
\end{equation}
\end{proposition}
\begin{corollary} If $\gamma \equiv \gamma^{\prime}$
then
$$(\alpha+\beta).\gamma \equiv \alpha.\gamma + \beta.\gamma^{\prime}.$$
\end{corollary}
\begin{proof}
Let $e^{\gamma t} \simeq h(t)$ the generating function of the umbra $\gamma.$
It is
$$e^{[(\alpha+\beta).\gamma]t} \simeq [h(t)]^{\alpha+\beta}
\simeq [h(t)]^{\alpha} [h(t)]^{\beta} \simeq
e^{(\alpha.\gamma)t} e^{(\beta.\gamma^{\prime})t}$$
from which the result follows. \qed
\end{proof}
\begin{remark}
{\rm As Taylor suggests in \cite{Taylor}, the auxiliary umbra
$\beta.\alpha$ provides an umbral interpretation of the random
sum since the moment generating function $g[\log f(t)]$
corresponds to the r.v. $S_N = X_1+X_2+\cdots+X_N$ with $X_i$
indipendent identically distributed (i.i.d.) r.v. having
moment generating function $f(t)$ and with $N$ a discrete
r.v. having moment generating function $g(t).$
The probabilistic interpretation of the corollary 2 states that the
random sum $S_{N+M}$ is similar to $S_{N} + S_{M},$ where $N$ and $M$
are two indipendent discrete r.v.
\par
The left distributive property of the point product respect
to the sum does not hold since
$$e^{[\alpha.(\beta+\gamma)]t} \simeq [g(t)]^{\alpha}
[h(t)]^{\alpha} \not \simeq f[\log g(t)] f[\log h(t)]$$
where $g(t) \simeq e^{\beta t}, h(t) \simeq e^{\gamma t}$
and $f(t) \simeq e^{\alpha t}.$
\par
Again this result runs in parallel with the probability theory.
In fact, let $Z=X+Y$ be a r.v. with $X$ and $Y$  two indipendent
r.v. The random sum $S_N = Z_1+Z_2+\cdots+Z_N,$ with
$Z_i$ i.i.d. r.v. similar to $Z,$ is not similar to the r.v. $S^{X}_N+
S^{Y}_N$ where $S^{X}_N=X_1+X_2+\cdots+X_N$ and $X_i$ i.i.d.
r.v. similar to $X,$ and where $S^{Y}_N=Y_1+Y_2+\cdots+Y_N$
and $Y_i$ i.i.d. r.v. similar to $Y.$}
\end{remark}
\begin{corollary}
\begin{description}
\item[\it{(i)}] If $\beta.\alpha \equiv \beta.\gamma$ then
$\alpha \equiv \gamma;$
\item[\it{(ii)}] if $c \in R$ then $\beta.(c \alpha) \equiv c (\beta.\alpha)$
for any two distinct umbrae $\alpha$ and $\beta;$
\item[\it{(iii)}] $\beta.(\gamma.\alpha) \equiv (\beta.\gamma).\alpha.$
\end{description}
\end{corollary}
\begin{proof}
Via generating functions. \qed
\end{proof}
To end this section, we deal with the notion of the inverse of an
umbra. Two umbrae $\alpha$ and $\beta$ are said to be {\it inverse}
to each other when $\alpha + \beta \equiv \varepsilon.$ Recall that,
dealing with a saturated umbral calculus, the inverse of an
umbra is not unique, but any two inverse umbrae of the umbra
$\alpha$ are similar.
\begin{proposition}
If $\alpha$ is an umbra with generating function $e^{\alpha t}
\simeq f(t)$ then its inverse $\beta$ has generating function
$e^{\beta t} \simeq [f(t)]^{-1}.$
\end{proposition}
\begin{proof}
The result follows observing that $e^{(\alpha+\beta)t} \simeq 1.$ \qed
\end{proof}
Similarly, for every positive integer $n$ and for every umbra $\alpha \in A,$
the inverse of the auxiliary umbra $n.\alpha,$ written as $-n.\alpha^{\prime}$
with $\alpha \equiv \alpha^{\prime},$ is similar to $\beta^{\prime} +
\beta^{\prime \prime} + \cdots + \beta^{\prime \prime \prime}$
where $\beta^{\prime}, \beta^{\prime \prime}, \ldots, \beta^{\prime \prime \prime}$
is any set of $n$ distinct umbrae similar to $\beta,$ being $\beta$
the inverse of $\alpha.$ The notation $-n.\alpha^{\prime}$ is
justified by noting that
$$n.\alpha-n.\alpha^{\prime} \equiv (n-n).\alpha \equiv 0.\alpha \equiv
\varepsilon.$$
\begin{proposition}
If $\alpha$ is an umbra with generating function $e^{\alpha t}
\simeq f(t),$ then the inverse of $n.\alpha$ has generating function
$e^{(-n.\alpha^{\prime}) t} \simeq [f(t)]^{-n}.$
\end{proposition}
\begin{proof}
The result follows observing that $e^{(n.\alpha-n.\alpha^{\prime})t} \simeq 1.$
\qed
\end{proof}
The inverse of the umbra $x.\alpha$ is the umbra $-x.\alpha^{\prime}$
with $\alpha \equiv \alpha^{\prime}$ defined by
$$
x.\alpha-x.\alpha^{\prime} \equiv (x-x).\alpha \equiv 0.\alpha \equiv
\varepsilon.$$
%
\subsection{The point power}
%
As it is easy to be expected, the definition of the power
of moments requires the use of similar umbrae and so of a point
operation. This notion comes into this picture by a
natural way, providing also an useful tool for umbral
manipulation of generating function.
\par
We shall denote by the symbol $\alpha^{.n}$
an auxiliary umbra similar to the product $\alpha^{\prime}
\alpha^{\prime \prime} \cdots \alpha^{\prime\prime\prime}$
where $\alpha^{\prime},\alpha^{\prime \prime},\ldots,\alpha^{\prime\prime
\prime}$ are a set of $n$ distinct umbrae each of which is similar
to the umbra $\alpha.$ We assume that $\alpha^{.0}$ is an umbra similar
to the unity umbra $u.$ A similar notion is introduced for any
umbral polynomial $p.$ The following statements are easily to
be proved:
\begin{proposition}
\smallskip
\begin{description}
\item[\it{(i)}] If $c \in R$ then $(c \alpha)^{.n} \equiv c^n \alpha^{.n}$
for any nonnegative integer $n \ne 0;$
\item[{\it (ii)}] $(\alpha^{.n})^{.m} \equiv \alpha^{.nm}
\equiv (\alpha^{.m})^{.n}$ for any two nonnegative integers $n,m;$
\item[{\it (iii)}]  $\alpha^{.(n+m)} \equiv \alpha^{.n}
(\alpha^{\prime})^{.m}$ for any two nonnegative integers $n,m$ and any two distinct umbrae
$\alpha \equiv \alpha^{\prime};$
\item[{\it iv)}] $(\alpha^{.n})^k \equiv (\alpha^k)^{.n}$ for any two nonnegative integers
$n,k.$
\end{description}
\end{proposition}
By the last statement, the moments of $\alpha^{.n}$ for any integer $n$ are:
\begin{equation}
E[(\alpha^{.n})^k] =E[(\alpha^{k})^{.n}] = a_k^n, \quad
k=0,1,2,\ldots
\label{(eq:10)}
\end{equation}
so that the moments of the umbra $\alpha^{.n}$ are the $n-$th
power of the moments of the umbra $\alpha.$
\begin{proposition}
The generating function of the $n-$th point power of the umbra $\alpha$
is the $n-$th power of the generating function of the umbra $\alpha.$
\end{proposition}
Note that, by virtue of Propositions 2 and 8 it is
\begin{equation}
e^{(n.\alpha)t} \simeq (e^{\alpha t})^{.n}.
\label{(eq:11)}
\end{equation}
The relation (\ref{(eq:11)}) restores the natural umbral interpretation
of $[f(t)]^n.$ More, let us observe that if $\alpha$ and $\beta$ are
not similar, it is
$$(\alpha+\beta)^{.n} \equiv \sum_{i=0}^n
\left( \begin{array}{c}
n \\
i
\end{array} \right) \alpha^{.i} \beta^{.(n-i)}.$$
\par
The point power operation leads us to define the point
exponential of an umbra. We shall denote by the symbol $e.^{\alpha} $
the auxiliary umbra
\begin{equation}
e.^{\alpha} \equiv u + \sum_{n=1}^{\infty}
\frac{\alpha^{.n}}{n!}.
\label{(eq:12)}
\end{equation}
We have immediately $e.^{\epsilon} \equiv u.$
\begin{proposition} For any umbra $\alpha,$ it is
\begin{equation}
e.^{(n.\alpha)} \simeq (e.^{\alpha})^{.n}.
\label{(eq:13)}
\end{equation}
\end{proposition}
\begin{proof}
It results
$$E[(e.^{\alpha})^{.n}]=E[e.^{\alpha}]^n = e^{n E[\alpha]}
= \sum_{k \geq 0} \frac{n^k E[\alpha]^k}{k!}$$
and also
$$E[e.^{(n.\alpha)}]=  \sum_{k \geq 0} \frac{E[(n.\alpha)^{.k}]}{k!}
=  \sum_{k \geq 0} \frac{n^k E[\alpha.^k]}{k!},$$
by which (\ref{(eq:13)}) follows. \qed
\end{proof}
\noindent \indent
Up to similarity, the expression of the moments given in (\ref{(eq:10)}) justifies the
definition of the auxiliary umbra $\alpha^{.x}$ as the umbra whose moments
are
$$E[(\alpha^{.x})^k]=a_k^x, \quad k=0,1,2,\ldots.$$
\begin{proposition}
Let $\alpha$ be an umbra and $e^{\alpha t}\simeq f(t)$ its
generating function. It is
$$e^{(x.\alpha)t} \simeq (e^{\alpha t})^{.x} \simeq [f(t)]^x.$$
\end{proposition}
Via moments, it is possible to prove the analogue of
Proposition 7 where $n$ and $m$ are replaced by $x$ and $y$
with $x,y \in R.$
\par
Once again, we define the auxiliary umbra $\alpha^{.\beta}$ as the umbra
whose moments are umbrally equivalent to
$$(\alpha^{.\beta})^k \simeq a_k^{\beta}, \quad k=0,1,2,\ldots$$
and we set $\epsilon^{.\alpha} \equiv \epsilon.$
\begin{proposition}
\begin{enumerate}
\item[{\it (i)}] $(\alpha^{.\beta})^{.\gamma} \equiv \alpha^{.(\gamma.\beta)};$
\item[{\it (ii)}] $\alpha^{.(\beta+\gamma)} \equiv \alpha^{.\beta}
(\alpha^{\prime})^{.\gamma}$ for any two distinct umbrae
$\alpha \equiv \alpha^{\prime}.$
\end{enumerate}
\end{proposition}
\begin{proof}
It follows via moments. \qed
\end{proof}
\begin{proposition}
Let $e^{\alpha t}\simeq f(t)$ be the generating function of the umbra $\alpha.$
It is
$$e^{(\beta.\alpha)t} \simeq (e^{\alpha t})^{.\beta} \simeq [f(t)]^\beta.$$
\end{proposition}
In closing, we notice that the generating function of the point
product between umbrae is umbrally equivalent to the following
series:
\begin{equation}
e^{(\beta.\alpha)t} \simeq \sum_{i=0}^{\infty} (\beta)_i
\frac{[e^{\alpha t}-u]^{.i}}{i!}
\label{(eq:14)}
\end{equation}
by the relation (\ref{(eq:8)}) and Proposition 8.
%
\section{Bell umbrae}
The Bell numbers $B_n$ have a long history and their origin
is unknown: Bell ascribes them to Euler
even without a specific reference \cite{Bell3}. Usually they are
referred as the number of the partitions of a finite nonempty set
with $n$ elements or as the coefficients of the Taylor series
expansion of the function $\exp(e^t-1).$ It is just writing about the Bell
numbers  that Gian-Carlo Rota \cite{Bellnumbers} gives the first glimmering
of the effectiveness of the umbral calculus in manipulating
number sequences, indeed his proof of the Dobinski's formula
is implicitly of umbral nature.
\par
In this section, the umbral definition of the Bell numbers
allows the proofs of several classical identities (cf.
\cite{Touchard}) through elementary arguments and smooths the way to the umbral
interpretation of the Poisson random variables.
\begin{definition}
An umbra $\beta$ is said to be a {\it Bell scalar umbra} if
$$(\beta)_n \simeq 1 \qquad n=0,1,2,\ldots$$
where $(\beta)_0=1$ and $(\beta)_n=\beta(\beta-1)\cdots(\beta-n+1)$
is the lower factorial.
\end{definition}
Up to similarity, the Bell number sequence is umbrally represented by the
Bell scalar umbra. Indeed, being
$$\beta^n=\sum_{k=0}^{n} S(n,k) (\beta)_k$$
where $S(n,k)$ are the Stirling numbers of second kind, then
$$E(\beta^n)=\sum_{k=0}^{n} S(n,k) E[(\beta)_k] =
\sum_{k=0}^{n} S(n,k) = B_n$$
where $B_n$ are the Bell numbers.
\par
The following theorem provides a characterization
of the Bell umbra.
\begin{theorem}
A scalar umbra $\beta$ is a Bell umbra iff
\begin{equation}
\beta^{n+1} \simeq (\beta+u)^n \quad n=0,1,2,\ldots
\label{(eq:16)}
\end{equation}
\end{theorem}
\begin{proof}
If $\beta$ is the Bell scalar umbra, being
$\beta(\beta-u)_n \simeq (\beta)_{n+1}$ it is $E[\beta(\beta-u)_n]=
1 = E[(\beta)_{n}].$ By the linearity it follows
$$E[\beta p(\beta-u)]=E[p(\beta)]$$
for every polynomial $p$ in $\beta.$ So, the idenity
(\ref{(eq:16)}) follows setting $p(\beta)=(\beta+u)^n.$
Viceversa, the relation (\ref{(eq:16)}) gives
$$E[\beta^{n+1}]= E[(\beta+u)^n]= \sum_{k=0}^n
\left( \begin{array}{c}
n \\
k
\end{array} \right) E[\beta^k]$$
or setting $E[\beta^n]=B_n$ one has
$$B_{n+1} =  \sum_{k=0}^n
\left( \begin{array}{c}
n \\
k
\end{array} \right)B_k$$
that is the recursion formula of the Bell numbers. \qed
\end{proof}
\begin{corollary}
If $\beta$ is the Bell scalar umbra, then
\begin{equation}
D_t[e^{\beta t}] \simeq e^{(\beta + u)t}.
\label{(eq:17)}
\end{equation}
\end{corollary}
\begin{proposition}
If $\beta$ is the Bell scalar umbra, for any integer $k>0$ and for any polynomial $p(x)$
the following relation holds
$$p(\beta+k.u) \simeq (\beta)_k p(\beta) \simeq p(\beta).$$
\end{proposition}
\begin{proof}
For $n \geq k,$ by the definition 1 it follows
$$(\beta)_n \simeq (\beta)_{n+k} \simeq (\beta)_k (\beta-k.u)_n.$$
Thus for any polynomial $q$ it is
$$q(\beta) \simeq (\beta)_k q(\beta-k.u)$$
by which one has
$$(\beta+k.u)^n \simeq  (\beta)_k \beta^{n}, \quad n=0,1,2,\ldots$$
setting $q(\beta)=(\beta+k.u)^n.$ The result follows by
linearity. \qed
\end{proof}
\begin{proposition}
The generating function of the Bell umbra is
\begin{equation}
e^{\beta t} \simeq e.^{e^{u t} - u}.
\label{(eq:18)}
\end{equation}
\end{proposition}
\begin{proof}
By the definition 1 and the relation (\ref{(eq:14)})
it is
$$e^{\beta t} \simeq e^{(\beta.u)t} \simeq \sum_{i=0}^{\infty}
\frac{[e^{u t}-u]^{.i}}{i!}.$$
Thus (\ref{(eq:18)}) follows from the relation (\ref{(eq:12)}).
\qed
\end{proof}
\begin{remark}
{\rm
Let us go on with our probabilistic counterpoint noting that
the Bell umbra can be view as a Poisson r.v. with parameter
$\lambda=1.$ Indeed,  the moment generating function of the
Bell umbra is $\exp(e^t-1)$ (see (\ref{(eq:18)})) so that $P(e^t) = \exp(e^t-1)$
where $P(t)$ is the probability generating function, and
therefore $P(s)=\exp(s-1).$ By that, the moments
of a Poisson r.v. with parameter $1$ are the Bell numbers and
its factorial moments are equal to $1.$}
\end{remark}
The following theorem makes clear how the proof of Dobinski's formula
becomes natural through the umbral expression of generating function.
\begin{theorem}[Umbral Dobinski's formula] The Bell umbra $\beta$
satisfies the following formula:
$$\beta^{n} \simeq e.^{-u} \sum_{k=0}^{\infty} \frac{(k.u)^n}{k!}.$$
\end{theorem}
\begin{proof}
Being $e^{\beta t} \simeq e.^{-u} e.^{e^{u t}}$ it is
$$e^{\beta t} \simeq e.^{-u} \sum_{k=0}^{\infty} \frac{e^{(k.u) t}}{k!}
\simeq e.^{-u} \sum_{k=0}^{\infty} \frac{1}{k!} \left\{ \sum_{n=0}^{\infty}
\frac{(k.u)^n t^n}{n!} \right\}$$
by which the result follows. \qed
\end{proof}
%
\subsection{The Bell polynomial umbra}
%
\begin{definition}
An umbra $\phi$ is said to be a {\it Bell polynomial umbra} if
$$(\phi)_n \simeq x^n \qquad n=0,1,2,\ldots.$$
\end{definition}
Note that $\phi \equiv \beta$ for $x=1.$ Moreover, being
$$\phi^n=\sum_{k=0}^{n} S(n,k) (\phi)_k,$$
by the definition 2 it follows
\begin{equation}
E(\phi^n)=\sum_{k=0}^{n} S(n,k) E[(\phi)_k] = \sum_{k=0}^{n} S(n,k) x^k =
\Phi_n(x).
\label{(eq:20)}
\end{equation}
The polynomials $\Phi_n(x)$ have a statistical origin and are known in the
literature as {\it exponential polynomials}. Indeed, they
were first introduced by Steffensen \cite{Steffensen} and studied further by
Touchard \cite{Touchard} and others. Rota, Kahaner and Odlyzko \cite{VIII}
state their basic properties via umbral operators.
\begin{proposition}
The generating function of the Bell polynomial umbra is
\begin{equation}
e^{\Phi t}  \simeq  e.^{x.(e^{u t} - u)}
\label{(eq:21)}
\end{equation}
\end{proposition}
\begin{proof}
By the definition 2 and the relation (\ref{(eq:14)})
it is
$$e^{\phi t} \simeq e^{(\phi.u)t} \simeq \sum_{i=0}^{\infty}
x^i \frac{[e^{u t}-u]^{.i}}{i!}.$$
Thus (\ref{(eq:21)}) follows from the relation (\ref{(eq:12)}).
\qed
\end{proof}
The following theorem provides a characterization of the Bell polynomial
umbra.
\begin{theorem}
An umbra $\phi$ is the Bell polynomial umbra iff
$$\phi \equiv x.\beta$$
where $\beta$ is the Bell scalar umbra.
\end{theorem}
\begin{proof}
The result comes via (\ref{(eq:21)}). \qed
\end{proof}
\begin{remark}
{\rm
The Bell polynomial umbra can be view as a Poisson r.v. with parameter
$\lambda=x.$ Indeed,  the moment generating function of the Bell polynomial umbra is
$\exp[x(e^t-1)]$ (see (\ref{(eq:21)})) so that $P(e^t) = \exp[x(e^t-1)],$
where $P(t)$ is the probability generating function,
and therefore $P(s)=\exp[x(s-1)].$ By that, the moments of a Poisson r.v.
with parameter $x$ are the exponential polynomials and its factorial
moments are equal to $x^n.$ \par When $x=n,$ the Bell polynomial
umbra $n.\beta$ is the sum of $n$ similar
uncorrelated Bell scalar umbrae, likewise in probability theory where a Poisson r.v. of
parameter $n$ can be view as the sum of $n$ i.i.d. (eventually uncorrelated)
Poisson r.v. with parameter $1.$ More in general, the closure under
convolution of the Poisson probability distributions i.e. $F_s \star F_t = F_{s+t},$ where
$F_t$ is a Poisson probability distribution depending on the parameter $t,$
is umbrally translated by $x.\beta+y.\beta^{\prime} \equiv (x+y).\beta$
(cf. statement {\it (iv)} of Proposition 1).}
\end{remark}
The next theorem is the polynomial analogue of the
theorem.
\begin{theorem}
An auxiliary umbra $x.\beta$ is a Bell polynomial umbra iff
\begin{equation}
(x.\beta)^{n+1} \simeq x (x.\beta+u)^n, \quad n=0,1,2,\ldots
\label{(eq:22)}
\end{equation}
\end{theorem}
\begin{proof}
Observe that $D_t[e^{(x.\beta)t}] \simeq D_t[ (e^{\beta t})^{.x}]
\simeq x e^{[(x-1).\beta]t} D_t[e^{\beta^{\prime} t}],$ where
$\beta^{\prime} \equiv \beta.$ From
(\ref{(eq:17)}) it is $D_t[e^{(x.\beta)t}] \simeq x e^{(x.\beta+u)t}$
from which the result follows immediately.
Viceversa, the relation (\ref{(eq:22)}) gives (\ref{(eq:16)})
for $x=1,$ by which it follows that $\beta$ is the Bell scalar
umbra. \qed
\end{proof}
The formula (\ref{(eq:22)}) represents the umbral
equivalent of the well known recursive formula for the exponential
polynomials:
$$\Phi_{n+1}(x)= x \sum_{k=0}^n \left( \begin{array}{c}
n \\
k
\end{array} \right)\Phi_k(x).$$
Similarly, the next proposition gives an umbral analogue of the Rodrigues formula
for the exponential polynomials (cf. \cite{VIII}).
\begin{proposition}
The Bell polynomial umbra $x.\beta$ has the following property:
$$D_x[(x.\beta)^n] \simeq (x.\beta + u)^n - (x.\beta)^n.$$
\end{proposition}
\begin{proof}
>From (\ref{(eq:18)}) it is
$$D_x[e^{(x.\beta)t}] \simeq e^{(x.\beta)t} (e^{ut}-u) \simeq e^{(x.\beta+u)t} -
e^{(x.\beta)t}$$
by which the result follows immediately. \qed
\end{proof}
In closing we state the polynomial version of the umbral Dobinski's formula.
\begin{proposition} The Bell polynomial umbra $x.\beta$ satisfies the following relation:
$$(x.\beta)^n \simeq e^{-x.u} \sum_{k=0}^{\infty} \frac{(k.u)^n x^k}{k!}.$$
\end{proposition}
\begin{proof}
Being $e^{(x.\beta) t} \simeq e.^{-x.u} e.^{(x.e^{u t})}$ it is
$$e^{(x.\beta) t} \simeq e.^{-x.u} \sum_{k=0}^{\infty} \frac{(x.e^{ut})^{.k}}{k!}
\simeq e.^{-x.u} \sum_{k=0}^{\infty} \frac{x^k e^{(k.u) t}}{k!}$$
by which the result follows. \qed
\end{proof}
%
\subsection{The exponential umbral polynomials}
%
Let us introduce a new family of umbral polynomials
that turns out to be an useful tool in the umbral composition, also
disclosing an unexpected probabilistic interpretation.
\par
Set
\begin{equation}
\Phi_n(\alpha) = \sum_{k=0}^n S(n,k) \alpha^k, \quad
n=0,1,2,\ldots,
\label{(eq:22-1)}
\end{equation}
we will call $\Phi_n(\alpha)$ {\it exponential umbral
polynomials}. By the identity (\ref{(eq:20)}), being
$(x.\beta)^n \simeq \Phi_n(x),$ where $\beta$ is the Bell scalar umbra, it
is glaring that
$$\Phi_n(\alpha) \simeq  (\alpha.\beta)^n \quad n=0,1,2,\ldots$$
and
\begin{equation}
(\alpha.\beta)_n \simeq \alpha^n \quad n=0,1,2,\ldots,
\label{(eq:22-2)}
\end{equation}
a formal proof passing through similar arguments already
produced for the umbra $x.\beta.$
\begin{proposition}
Let $\beta$ be the Bell scalar umbra. If $e^{\alpha t}
\simeq f(t)$ is the generating function of the umbra $\alpha$ then
\begin{equation}
e^{(\alpha.\beta) t}  \simeq  f[e^t-1].
\label{(eq:22-3)}
\end{equation}
\end{proposition}
\begin{proof}
The result follows by the relation (\ref{(eq:9)}) observing that
$e^{\beta t} \simeq e^{e^t-1}.$ \qed
\end{proof}
When $f(t)$ is considered as the moment generating function of a r.v. $X,$
a probabilistic interpretation of (\ref{(eq:22-3)}) suggests that the umbra
$\alpha.\beta$ represents a Poisson r.v. $N_X$ with random
parameter $X.$  Indeed the probability generating function of $N_X$ is
\begin{eqnarray*}
P(s) & = & \sum_{k=0}^{\infty} {\bf P}(N_X=k) s^k = \sum_{k=0}^{\infty}
s^k \int_0^{\infty} {\bf P}(N_X=k| X=x) dF_X(x) \\
& = & \sum_{k=0}^{\infty} \frac{s^k}{k!} \int_0^{\infty} x^k e^{-x} dF_X(x)
= \sum_{k=0}^{\infty} \frac{(s-1)^k}{k!} E[X^k] = f(s-1)
\end{eqnarray*}
hence the moment generating function of $N_X$ is $f(e^t-1).$
To the best of our knowledge, this r.v. has been introduced in
\cite{FellerII} as the randomized Poisson r.v. Once more
the closure under convolution of the Poisson probability distributions
leads us to claim that the point product $\alpha.\beta$ is the
umbral corresponding of the random sum of independent Poisson r.v. with
parameter $1$ indexed by an integer r.v. $X.$
%
\section{The partition umbra}
%
As suggested in \cite{VIII}, there is a connection between
polynomials of binomial type and compound Poisson processes.
Two different approaches can be found in \cite{Stam1} and in
\cite{Cerasoli}. In this section, we suggest a way, that we believe
to be natural, in order to make clear this connection.
\begin{definition}
An umbra $\psi$ is said to be an {\it $\alpha-$partition umbra} if
$$\psi \equiv \beta.\alpha$$
with $\beta$ the Bell scalar umbra.
\end{definition}
Note that the $u-$partition umbra is the Bell scalar umbra.
\begin{proposition}
The generating function of the $\alpha-$partition umbra $\psi$ is
\begin{equation}
e^{\psi t} \simeq e.^{(e^{\alpha t}  - u)}.
\label{(eq:24)}
\end{equation}
\end{proposition}
\begin{proof}
>From (\ref{(eq:14)}) and by the definition 1 one has
$$e^{(\beta.\alpha)t}  \simeq
 \sum_{n=0}^{\infty} \frac{(e^{\alpha t}-u)^{.n}}{n!}. $$
Thus (\ref{(eq:24)}) follows from the relation (\ref{(eq:12)}).
\qed
\end{proof}
\noindent \indent
The generating function (\ref{(eq:24)}) leads us to interpret a partition umbra
as a compound Poisson r.v. with parameter $1.$ As well known (cf. \cite{FellerI}), a
compound Poisson r.v. with parameter $x$ is  introduced as a random sum $S_N = X_1 + X_2 +
\cdots + X_N$ where $N$ has a Poisson distribution with parameter $x$. The point product of
definition 3 fits perfectly this probabilistic notion
taking into consideration that the Bell scalar umbra $\beta$ plays the role of
a Poisson r.v. with parameter $1.$ What's more, since the Poisson r.v. with parameter
$x$ is umbrally represented by the Bell polynomial umbra $x.\beta,$
a compound Poisson r.v. with parameter $x$ is represented by the
{\it polynomial $\alpha-$partition umbra} $x.\psi \equiv x.\beta.\alpha$
with generating function
\begin{equation}
e^{(x.\psi)t} \simeq e.^{[x.(e^{\alpha t}  - u)]}.
\label{(eq:25)}
\end{equation}
The name \lq \lq partition umbra\rq\rq$\,$ has also a probabilistic ground.
Indeed the parameter of a Poisson r.v. is usually denoted by $x=\lambda t,$
with $t$ representing a time interval, so that when this interval is
partitioned into non-overlapping ones, their contributions
are stochastic independent and add to $S_N.$
The last circumstance is umbrally expressed by the relation
\begin{equation}
(x+y).\beta.\alpha \equiv x.\beta.\alpha + y.\beta.\alpha
\label{(eq:somma)}
\end{equation}
that also assures the binomial property for the polynomial sequence
defined by $x.\beta.\alpha.$  In terms of generating functions, the
formula (\ref{(eq:somma)}) means that
\begin{equation}
h_{x+y}(t) = h_x(t)h_y(t)
\label{(eq:somma1)}
\end{equation}
where $h_x(t)$ is the generating function of $x.\beta.\alpha.$
Viceversa every generating function $h_x(t)$ satisfying the
equality (\ref{(eq:somma1)}) is the generating function of
a polynomial $\alpha-$partition umbra, namely $h_x(t)$ has an umbral
expression of the form (\ref{(eq:25)}).
\par
Going back the moments of a partition umbra, according to the
definition of the Bell scalar umbra and from (\ref{(eq:7)}) it
is
\begin{equation}
E[(\beta.\alpha)^n] = \sum_{k=1}^n B_{n,k}(a_1,a_2,\ldots,
a_{n-k+1})= Y_n(a_1,a_2,\ldots,a_n)
\label{(eq:26)}
\end{equation}
where $Y_n = Y_n(a_1,a_2,\ldots,a_n)$ are the {\it partition
polynomials} (or complete Bell exponential polynomials) and $a_i$ are
the moments of the umbra $\alpha.$ Although the complexity of the
partition polynomial expression, their umbral interpretation
(\ref{(eq:26)}) allows an easy proof that they are of binomial
type, simply observing that $\beta.\alpha + \beta.\gamma \equiv
\beta.(\alpha + \gamma).$
\par
Partition polynomials have been first introduced  by Bell \cite{Bell0}
who gave a pioneer umbral version of them in \cite{Bell1}.
Because of their generality, they include a variety of other
polynomials such as the cycle indicator of the symmetric
group and other of interest in number theory.
\par
As already done for the Bell scalar umbra, the next
theorem characterizes the partition umbrae and also provides
the following recursive formula for the partition
polynomials:
$$Y_{n+1}(a_1,a_2,\cdots,a_{n+1}) = \sum_{k=0}^n \left( \begin{array}{c}
n \\
k
\end{array} \right) a_{n-k+1} Y_k(a_1,a_2,\cdots,a_k).$$
\begin{theorem}
Every $\alpha-$partition umbra verifies the following relation
\begin{equation}
(\beta.\alpha)^{n+1} \simeq \alpha^{\prime}(\beta.\alpha +
\alpha^{\prime})^n \quad \alpha^{\prime} \equiv \alpha, \,\, n=0,1,2,\ldots
\label{(eq:27)}
\end{equation}
and viceversa.
\end{theorem}
\begin{proof}
Let $\psi$ an $\alpha-$partition umbra. Then from (\ref{(eq:24)}) it is
$D_t[e^{(\beta.\alpha)t}] \simeq e^{(\beta.\alpha)t} D_t[e^{\alpha^{\prime} t}],$
where $\alpha^{\prime} \equiv \alpha.$ The identity (\ref{(eq:27)}) follows
observing that $D_t[e^{\alpha^{\prime} t}] \simeq \alpha^{\prime} e^{\alpha^{\prime} t}.$
Going back the previous steps, from (\ref{(eq:27)}) one has
that $\beta.\alpha$ has generating function (\ref{(eq:24)}) and
so it is an $\alpha-$partition umbra. \qed
\end{proof}
The moments of the polynomial $\alpha-$partition umbra are
\begin{eqnarray}
E[(x.\beta.\alpha)^n] & = & \sum_{k=1}^n (x.\beta)_k B_{n,k}(a_1,a_2,\ldots,
a_{n-k+1}) \nonumber \\
& = & \sum_{k=1}^n x^k B_{n,k}(a_1,a_2,\ldots, a_{n-k+1})
\label{(eq:28)}
\end{eqnarray}
according to the definition 2. The same arguments given in the proof
of the theorem 6 lead to state that every polynomial partition umbra
verifies the following formula
$$(x.\beta.\alpha)^{n+1} \simeq x\alpha^{\prime}(x.\beta.\alpha +
\alpha^{\prime})^n \quad \alpha^{\prime} \equiv \alpha, \,\, n=0,1,2,\ldots$$
and viceversa.
%
\subsection{Umbral expression of the functional composition}
%
An umbral wording of the functional composition of exponential formal power series
is a thorny matter. It was broached by Rota, Shen and Taylor in \cite{RotaTaylor}
passing through the sequence of Abel polynomials. In this last section,
we give an intrinsic umbral expression of this operation via the notion
of partition umbra.
\begin{definition}
A {\it composition umbra} of the umbrae $\alpha$ and $\gamma$
is the umbra
$$\chi \equiv \gamma.\beta.\alpha$$
where $\beta$ is the Bell scalar umbra.
\end{definition}
In other words, the composition umbra $\chi$ is the point
product of the umbra $\gamma$ and the $\alpha-$partition umbra
$\beta.\alpha.$
\begin{remark}
{\rm As already stressed in section 3.2, the umbra $\gamma.\beta$
represents a randomized Poisson r.v.. Hence it is natural to
look at the composition umbra as a new r.v. that we will call
{\it compound randomized Poisson r.v.} Moreover, being $(\gamma.\beta).\alpha
\equiv \gamma.(\beta.\alpha)$ (cf. statement {\it (ii)} of corollary
3), the previous relation allows to see this new r.v. from another
side: the umbra $\gamma.(\beta.\alpha)$ generalizes the concept of a random
sum of i.i.d. compound Poisson r.v. with parameter $1$ indexed
by an integer r.v. $X,$ i.e. a randomized compound Poisson r.v.
with random parameter $X.$}
\end{remark}
\begin{proposition}
The generating function of the composition umbra $\gamma.\beta.\alpha$
is the functional composition of the generating functions $e^{\alpha t} \simeq f(t)$
and $e^{\gamma t} \simeq g(t).$
\end{proposition}
\begin{proof}
Via (\ref{(eq:24)}) it is $e^{(\beta.\alpha) t} \simeq e^{f(t)-1}.$ The result
follows by (\ref{(eq:9)}) observing that $e^{[\gamma.(\beta.\alpha)] t} \simeq
g \left\{\log[e^{f(t)-1}]\right\}.$ \qed
\end{proof}
The moments of the composition umbra are
\begin{equation}
(\gamma.\beta.\alpha)^n \simeq \sum_{k=0}^n \gamma^k B_{n,k}(a_1,a_2,\ldots,a_{n-k+1})
\label{(eq:30)}
\end{equation}
where $a_i$ are the moments of the umbra $\alpha.$ Indeed, by
(\ref{(eq:7)}) it is
$$(\gamma.\beta.\alpha)^n \simeq \sum_{k=0}^n (\gamma.\beta)_k
B_{n,k}(a_1,a_2,\ldots,a_{n-k+1})$$
and (\ref{(eq:30)}) follows from (\ref{(eq:22-2)}).
\par
Once more, we give a characterization of the composition umbra in the next
theorem.
\begin{theorem}
Every composition umbra verifies the following relation
\begin{equation}
(\gamma.\beta.\alpha)^{n+1} \equiv \gamma \alpha^{\prime} (\gamma.\beta.\alpha +
\alpha^{\prime})^n \quad \alpha \equiv \alpha^{\prime}, \,\, n=0,1,2,\ldots
\label{(eq:31)}
\end{equation}
and viceversa.
\end{theorem}
\begin{proof}
Let $\chi$ a composition umbra of $\alpha$ and $\gamma$. Then
from Proposition 20, it is $D_t[e^{\chi t}] \simeq g^{\prime}[f(t)-1]f^{\prime}(t).$
Equation (\ref{(eq:31)}) follows being $f^{\prime}(t) \simeq \alpha^{\prime}
e^{\alpha^{\prime} t}$ with $\alpha^{\prime} \equiv \alpha$ and
$g^{\prime}[f(t)-1] \simeq \gamma e^{\chi t}.$
Going back the previous steps, from (\ref{(eq:31)}) it follows
that $\gamma.\beta.\alpha$ has generating function $g[f(t)-1]$ and
so it is a composition umbra of $\alpha$ and $\gamma.$ \qed
\end{proof}
At this point, as custom, we put to test the definition 4 of
composition umbra, giving a proof of Lagrange inversion
formula. In the literature (cf. \cite{dibucchianico} for a plenty of references)
different forms of the Lagrange inversion formula are derived
using umbral calculus. The main tool of our proof is the
umbral expression of the (partial) Bell exponential polynomials
that we state in the next proposition.
\begin{lemma}
It is
\begin{equation}
B_{n,k}(a_1,a_2,\ldots,a_{n-k+1}) \simeq \left( \begin{array}{c}
n \\
k
\end{array}\right) \alpha^{.k} (k.\overline{\alpha})^{n-k}
\label{(eq:32)}
\end{equation}
where $\overline{\alpha}$ is the umbra with moments $\displaystyle{E\left[\overline{\alpha}^{n} \right]=
\frac{a_{n+1}}{a_1 (n+1)}}, n=1,2,\ldots.$
\end{lemma}
\begin{proof}
By the identity (\ref{(eq:3)}) it results
$$B_{n,k}(a_1,a_2,\ldots,a_{n-k+1})  = \frac{1}{k!} D^{(n)}_t
[(f(t)-1)^{k}]_{t=0}$$
where $D^{(n)}_t[\cdot]_{t=0}$ is the $n-$th derivative respect
to $t$ evaluated in $t=0$ and $f(t) \simeq e^{\alpha t},$ so
that
$$B_{n,k}(a_1,a_2,\ldots,a_{n-k+1})  \simeq  \frac{1}{k!} D^{(n)}_t
[(e^{\alpha t} - u)^{.k}]_{t=0}.$$
On the other hand, by the moment expression of umbra $\overline{\alpha}$
it follows $e^{\alpha t} -u \simeq \alpha \, t \, e^{\overline{\alpha} t}.$
Therefore one has
\begin{eqnarray*}
D^{(n)}_t [(e^{\alpha t} - u)^{.k}] & \simeq & \alpha^{.k}
D^{(n)}_t[t^k \, e^{k.\overline{\alpha}}] \\
& \simeq &
\alpha^{.k} \sum_{j=0}^k \left( \begin{array}{c}
n \\
j
\end{array}\right) D^{(j)}_t [t^k] D^{(n-j)}_t [e^{k.\overline{\alpha}}],
\end{eqnarray*}
using the binomial property of the derivative operator.
Finally, the result follows evaluating the right hand side of
the previous formula in $t=0$ and observing that
$D^{(n-k)}_t[e^{k.\overline{\alpha}}]_{t=0} \simeq
(k.\overline{\alpha})^{n-k}.$ \qed
\end{proof}
\begin{remark}
{\rm Let $\alpha \equiv u.$ Then $\alpha^{.k} \simeq 1$ for $k=0,1,2,\ldots$
and $B_{n,k}(1,1,\ldots,1) = S(n,k)$ the Stirling number
of the second kind. Moreover it is $\overline{\alpha} \equiv (-1.\delta)$
where $\delta$ is the Bernoulli umbra whose moments are the Bernoulli numbers
(cf. \cite{RotaTaylor}). From Lemma 1 it results
$$S(n,k) \simeq \left( \begin{array}{c}
n \\
k
\end{array}\right)(-k.\delta)^{n-k}$$
as already stated by Rota and Taylor through a different approach
(cf. Proposition 9.1 \cite{SIAM})}
\end{remark}
\begin{theorem}[Lagrange inversion formula]
Let $e^{\alpha t} \simeq f(t)$ and $e^{\gamma t} \simeq g(t).$
If $g[f(t)-1]=f[g(t)-1]=1+t$ then
\begin{equation}
\alpha^{.k} \gamma^k \simeq (-k.\overline{\alpha})^{k-1},
k=1,2,\ldots.
\label{(eq:34)}
\end{equation}
\end{theorem}
\begin{proof}
By the formulas (\ref{(eq:30)}) and (\ref{(eq:32)}), it is
\begin{equation}
\chi^n \simeq \sum_{k=0}^n \left( \begin{array}{c}
n \\
k
\end{array}\right) \alpha^{.k} \gamma^k
(-k.\overline{\alpha})^{n-k}.
\label{(eq:35)}
\end{equation}
On the other hand, the Abel identity (\ref{(Abel)}) gives
\begin{equation}
\chi^n \simeq \sum_{k=0}^n \left( \begin{array}{c}
n \\
k
\end{array}\right) \chi (\chi - k.\overline{\alpha})^{k-1}
(k.\overline{\alpha})^{n-k}.
\label{(eq:36)}
\end{equation}
Comparing (\ref{(eq:35)}) with (\ref{(eq:36)}) one has
$$ \alpha^{.k} \gamma^k \simeq \chi (\chi -k.\overline{\alpha})^{k-1}$$
by which the result follows expanding the right hand side of
the previous formula by the binomial theorem and observing that
from $g[f(t)-1]=1+t$ it is $\chi \simeq 1$ and
$\chi^j \simeq 0, j=2,3,\ldots.$ \qed
\end{proof}
More explicitly, the formula (\ref{(eq:34)}) says that the $k-$th
coefficient of the generating function $g(t)$ is equal to the $(k-1)-$th coefficient
of the generating function $[(f(t)-1)/t]^{-k},$ when $g[f(t)-1]=1+t.$
Note that if $f(t)-1=t e^{-t}$ then $\alpha^{.k} \simeq 1, \,
\overline{\alpha} \equiv -1.u$ and from (\ref{(eq:34)}) it is $\gamma^k \simeq (k.u)^{k-1}.$
\par
In closing, let us observe that if $a_1=1$ then
$f(t)-1 \simeq t e^{\overline{\alpha} t}$ and the Lagrange
inversion formula (\ref{(eq:34)}) becomes
$$\gamma^{k} \equiv (-k.\overline{\alpha})^{k-1}.$$
On the other hand, if the generating function $g(t)$ is written as
$g(t)-1 \simeq t e^{\overline{\gamma} t}$ then
the Lagrange inversion formula (\ref{(eq:34)}) becomes
$$k \overline{\gamma}^{k-1} \simeq (-k.\overline{\alpha})^{k-1}$$
that is equivalent to the version given in \cite{RotaTaylor} by using the Abel
polynomial sequence and its delta operator.

\end{document}